\documentclass[11pt]{article}

 \usepackage{setspace}
\usepackage{natbib}

\usepackage{latexsym}
\usepackage{calc}

\usepackage{amssymb}
\usepackage{amsmath}
\usepackage{graphicx}

\usepackage{psfrag}

\usepackage{float}

\usepackage{color}

\usepackage{nicefrac}

\usepackage{ntheorem}

\usepackage{hyperref}

\usepackage{tikz}
\usetikzlibrary{backgrounds}

\newcounter{hours}\newcounter{minutes}

\def\nr{\par \noindent}

\def\beq{\begin{equation}}
\def\eeq{\end{equation}}

\newtheorem{theorem}{Theorem}

\newtheorem{lemma}{Lemma}
\newtheorem{corollary}{Corollary}

\newtheorem{proposition}{Proposition}
\newtheorem{assumption}{Assumption}
\newtheorem{definition}{Definition}

\newtheorem{example}{Example}
\newtheorem{remark}{Remark}
\newcommand{\proof}{\bf Proof: \rm \nr}

\newcommand{\qed}{\hfill $\Box$ \nr \medskip}

\def\ba{\begin{array}}
\def\ea{\end{array}}
\def\beann{\begin{eqnarray*}}
\def\eeann{\end{eqnarray*}}
\def\bea{\begin{eqnarray}}
\def\eea{\end{eqnarray}}

\def\BT{\begin{theorem}}
\def\ET{\end{theorem}}
\def\BL{\begin{lemma}}
\def\EL{\end{lemma}}
\def\BC{\begin{corollary}}
\def\EC{\end{corollary}}
\def\BE{\begin{example}}
\def\EE{\end{example}}
\def\BD{\begin{definition}}
\def\ED{\end{definition}}
\def\BR{\begin{remark}}
\def\ER{\end{remark}}
\def\BAS{\begin{assumption}}
\def\EAS{\end{assumption}}
\def\BI{\begin{itemize}}
\def\EI{\end{itemize}}

\def\BMP{\begin{minipage}{9.5cm}}
\def\EMP{\end{minipage}}
\def\MPT{\begin{minipage}{11.5cm}}
\def\EPT{\end{minipage}}

\def\R{\mathbb{R}}

\newcommand*\samethanks[1][\value{footnote}]{\footnotemark[#1]}
\newcommand{\NNorm}[2]{\left\Vert {#1} \right\Vert_{#2}}
\newcommand{\AV}[1]{\left\vert {#1} \right\vert}

\title{Critical point theory for sparse recovery}
\author{
S. L\"ammel
\thanks{
Department of Mathematics, Chemnitz University of Technology,
Reichenhainer Str. 41, 09126
Chemnitz, Germany; e-mail: sebastian.laemmel@mathematik.tu-chemnitz.de, vladimir.shikhman@mathematik.tu-chemnitz.de (corresponding author).
 } \and V. Shikhman\samethanks[1]
}

\begin{document}
\maketitle
\vspace{-5ex}
\abstract{ 
We study the problem of sparse recovery in the context of compressed sensing.
This is to minimize the sensing error of linear measurements by sparse vectors  with at most $s$ non-zero entries.
We develop the so-called critical point theory for sparse recovery. This is done by introducing nondegenerate M-stationary points which adequately describe the global structure of this nonconvex optimization problem. We show that all M-stationary points are generically nondegenerate. In particular, the sparsity constraint is active at all local minimizers of a generic sparse recovery problem. Additionally, the equivalence of strong stability and nondegeneracy for M-stationary points is shown. 
We claim that the appearance of saddle points - these are M-stationary points with exactly $s-1$ non-zero entries - cannot be neglected. For this purpose we derive a so-called Morse relation, which gives a lower bound on the number of saddle points in terms of the number of local  minimizers. The relatively involved structure of saddle points can be seen as a source of well-known difficulty by solving the problem of sparse recovery to global optimality. 
}

\vspace{2ex}

\section{Introduction}
\label{sec:intro}
Compressed sensing is concerned with the recovery of a sparse vector $x$ from linear measurements $Ax=b$, where $A \in \R^{m \times n}$ is a sensing matrix and $b \in \R^m$ is a measurement vector. For this purpose, it is usual to consider the following optimization problem, see e.\,g. \cite{kutyniok:2012}:
\begin{equation}
\label{eq:cs1}
\min_{x \in \R^n} \,\, \left\|x\right\|_0\quad \mbox{s.\,t.} \quad Ax=b,
\end{equation}
where the so-called $\ell_0$ "norm" counts non-zero entries of $x$, i.\,e.
\[
\NNorm{x}{0}= \AV{\left\{i\in \{1,\ldots,n\} \,\left\vert\, x_i \ne 0 \right.\right\}}.
\]
If the linear measurements are prone to Gaussian noise, the optimization problem (\ref{eq:cs1}) can be modified as follows:
\begin{equation}
\label{eq:cs2}
\min_{x \in \R^n} \,\, \left\|x\right\|_0\quad \mbox{s.\,t.} \quad \left\|Ax-b\right\|_2 \leq \varepsilon,
\end{equation}
where $\varepsilon > 0$ is the bound on the error magnitude
with respect to the Euclidean norm. 
In this paper we consider an analogue formulation of (\ref{eq:cs2}) introduced by \cite{beck:2013}:
\[
\mbox{SR}:\quad \min_{x \in \R^n} \,\, \frac{1}{2}\left\|Ax-b\right\|^2_2\quad \mbox{s.\,t.} \quad \NNorm{x}{0} \le s,
\]
where $s \in\{0,1, \ldots,n-1\}$ is the bound on the number of non-zero entries of $x$. Note that the {\it sparse recovery problem} consists of minimizing the sensing error by sparse vectors with at most $s$ non-zero entries. Sometimes we write $\mbox{SR}(A,b)$ for the problem of sparse recovery, in order to highlight the dependence on the data $(A,b)$. 
Throughout the paper we make the following not very restrictive assumption, cf. \cite{beck:2013}.

\begin{assumption}
\label{ass:cs}
The bound on the number of non-zero entries does not exceed the number of measurements in SR, i.\,e. $s \leq m$.
\end{assumption}

The difficulty of solving SR comes from the combinatorial nature of the sparsity constraint $\NNorm{x}{0} \le s$. 
Although the objective function $f$ of SR is convex, its feasible set is non-convex as a union of linear subspaces. Nevertheless, several attempts to tackle SR have been undertaken in recent years.   

In the seminal paper \cite{beck:2013} a generalization of SR with an arbitrary smooth objective function is considered. The latter is referred to by the authors as {\it sparsity constrained nonlinear optimization}:
\begin{equation}
\label{eq:scno}
  \min_{x \in \R^n} \,\, f(x)\quad \mbox{s.\,t.} \quad \NNorm{x}{0} \le s.
\end{equation}
The notions of basic feasibility, $L$-stationarity, and CW-minimality have been introduced and shown to be necessary optimality conditions for (\ref{eq:scno}). The formulation of $L$-stationarity mimics the techniques from convex optimization by using the orthogonal projection on the feasible set. The notion of CW-minimum incorporates the coordinate-wise optimality along the axes. Based on both stationarity concepts, algorithms that find points satisfying these conditions have been developed. These are the iterative hard thresholding method, as well as the greedy and partial sparse-simplex methods. In a series of subsequent papers \cite{beck:2016, beck:2018} elaborated the algorithmic approach based on $L$-stationarity and CW-minimality. 

Another line of research started with \cite{burdakov:2016}. Here, in addition to an arbitrary smooth objective function also smooth equality and inequality constraints have been incorporated into the feasible set. For that, the authors coin the new term of {\it mathematical programs with cardinality constraints}:
\begin{equation}
\label{eq:mpcc}
  \min_{x \in \R^n} \,\, f(x)\quad \mbox{s.\,t.} \quad \NNorm{x}{0} \le s, \quad h(x)=0, \quad g(x) \geq 0.
\end{equation}
The key idea in \cite{burdakov:2016} is to provide a mixed-integer formulation of (\ref{eq:mpcc}) whose standard relaxation still has the same solutions. For the relaxation the notion of S-stationary points is proposed. S-stationarity corresponds to the standard Karush-Kuhn-Tucker condition for the relaxed program. The techniques applied follow mainly those for mathematical programs with complementarity constraints. In particular, an appropriate regularization method for solving (\ref{eq:mpcc}) is suggested. The latter is proved to converge towards so-called M-stationary points. M-stationarity corresponds to the standard Karush-Kuhn-Tucker condition of the tightened program, where zero entries of a feasible point remain locally vanishing. Further research in this direction is presented in a series of subsequent papers \cite{cervinka:2016}, \cite{bucher:2018}.

The aim of this paper is to develop a {\it critical point theory} for the problem of sparse recovery. The main idea of critical point theory is to identify stationary points which roughly speaking induce the global structure of the underlying optimization problem. They have not only to include minimizers, but also all kinds of saddle points -- just in analogy to the unconstrained case. Critical point theory for other non-convex optimization problems, such as e.\,g. mathematical programs with complementarity constraints, general semi-infinite programming, mathematical problems with vanishing constraints, is elaborated in \cite{jongen:2009}, \cite{jongen:2011}, \cite{dorsch:2012}, respectively.

Let us overview our main results on the the critical point theory for SR:
\begin{itemize}
    \item[(i)] It turns out that the concept of M-stationarity from \cite{burdakov:2016} is the adequate stationarity concept for our purposes. We introduce the notion of {\it nondegeneracy} for M-stationary points of SR. It is proved that all M-stationary points are generically nondegenerate, see Theorem \ref{thm:generic-cs}. As an important consequence, the sparsity constraint must be active at all local minimizers of a generic SR, see Corollary \ref{cor:min}. 
    \item[(ii)] Further, we introduce the notion of {\it strongly stability} of M-stationary points in the sense of \cite{kojima:1980}. The equivalence of strong stability and nondegeneracy for M-stationary points of SR is shown, see Theorem \ref{thm:stab-nondeg}. In case of degeneracy a local minimizer of SR may bifurcate into multiple minimizers and a saddle point, see Example \ref{ex:cw-inst}.
    \item[(iii)] The role of saddle points play M-stationary points with exactly $s-1$ non-zero entries. We derive a so-called Morse relation, which gives a lower bound on the number of saddle points in terms of the number of local  minimizers, see Theorem \ref{thm:mrel}. Hence, the appearance of saddle points cannot be neglected at least from the perspective of global optimization. 
 As further novelty, a saddle point may lead to more than two different local minimizers. 
The relatively involved structure of saddle points can be seen as a source of well-known difficulty if solving mathematical programs with sparsity constraint to global optimality. 
\end{itemize}

We would like to mention that in the recent preprint \cite{laemmel:2019} the critical point theory for sparsity constrained nonlinear optimization (\ref{eq:scno}) has been established. 
Note that although SR constitutes a subclass of (\ref{eq:scno}), the adjustment of results from \cite{laemmel:2019} for SR is by far not straight-forward. In fact, we cannot use either the corresponding results or their proof technique, in order to show the genericity of nondegenerate M-stationary points of SR. This is due to the fact that the data space of $(A,b)$ generates just a subset of $C^2$-functions via 
\[
f(x)=\frac{1}{2}\left\|Ax-b\right\|^2_2.
\]
The issue of strong stability is new and has not been studied in \cite{laemmel:2019}. So is its equivalence to nondegeneracy for M-stationary points of SR. Finally, the derivation of Morse relation needs an SR specific notion of $s$-regularity of the sensing matrix $A$ being introduced by \cite{beck:2013}, see Lemma \ref{lem:help}. For the readers' convenience, we decided to make the exposition of the critical point theory for SR self-contained. This allows a potential reader, which is just interested in the topic of sparse recovery, not to consult the previous paper at all.

The paper is organized as follows. In Section \ref{sec:m-nond} we discuss the notion of a nondegenerate M-stationary point. In Section \ref{sec:gen} we show that nondegeneracy is a generic property of M-stationary points. Section \ref{sec:stab} is devoted to the strong stability of M-stationary points and its equivalence to their nondegeneracy. The global structure of SR is described in Section \ref{sec:glob}.

Our notation is standard. The cardinality of a finite set $S$ is denoted by $|S|$. The $n$-dimensional Euclidean space is denoted by $\R^n$ with the coordinate vectors $e_i$, $i=1, \ldots,n$. For $J \subset \{1, \ldots, n\}$ we denote by $\mbox{conv}\left(e_j, j \in J\right)$ the convex hull of the coordinate vectors $e_j, j \in J$. 
The Euclidean norm of a vector $x\in \R^n$ is denoted by 
$\|x\|_2$, and by $x \geq 0$ we refer to the componentwise comparison $x_i \geq 0$ for all $i=1,\ldots, n$.
The entries of the subvector $x_I$ correspond to those of $x \in \R^n$ with respect to a given index set $I\subset\{1, \ldots, n\}$.
The space of real $(m,n)$-matrices is denoted by $\R^{m\times n}$. For $A \in \R^{m\times n}$ the transposed matrix is denoted by $A^T \in \R^{n\times m}$. If $A \in \R^{m\times n}$ is of full rank $n \leq m$, then $A^+=(A^T A)^{-1}A^T$ denotes the Moore-Penrose inverse of $A$.  
For an index set $I \subset \{1, \ldots,n\}$ we denote by $A_I$ the submatrix of $A \in \R^{m\times n}$ with the columns corresponding to the set $I$. Additionally, we 
denote by $A_I^T$ the transposition of $A_I$.

\section{Nondegeneracy}
\label{sec:m-nond}
For $0 \leq k \leq n$ we use the notation
\[
  \R^{n,k}= \left\{ x \in \R^n\, \left\vert \, \NNorm{x}{0} \le k \right.\right\}.
\]
Using the latter, the feasible set of SR can be written as
\[
  \R^{n,s}= \left\{ x \in \R^n\, \left\vert \, \NNorm{x}{0} \le s \right.\right\}.
\]
For a feasible point $x\in \R^{n,s}$ we define the following complementary index sets:
\[
I_0(x) = \left\{i\in \{1,\ldots,n\} \,\left\vert\, x_i = 0 \right.\right\}, \quad
  I_1(x) = \left\{j\in \{1,\ldots,n\} \,\left\vert\, x_j \ne 0 \right.\right\}.
\]
Without loss of generality, we assume throughout the whole paper that
at the particular point of interest $\bar x \in  \R^{n,s}$ it holds:
\[
  I_0\left(\bar x\right) = \left\{1,\ldots,n-\NNorm{\bar x}{0}\right\}, \quad
  I_1\left(\bar x\right) = \left\{n-\NNorm{\bar x}{0}+1,\ldots,n\right\}.
\]
Using this convention, the following local description of SR feasible set can be deduced. Let $\bar x \in \R^{n,s}$ be a feasible point of SR. Then, there exist neighborhoods $U_{\bar x}$ and $V_0$ of $\bar x$ and $0$, respectively, such that under the linear coordinate transformation $\Phi(x)=x-\bar x$ we have:
\begin{equation}
\label{eq:diff}
\Phi\left(\R^{n,s}\cap U_{\bar x}\right) = \left( \R^{n-\NNorm{\bar x}{0},s-\NNorm{\bar x}{0}}\times \R^{\NNorm{\bar x}{0}} \right) \cap V_0, \quad \Phi\left(\bar x\right)=0.
\end{equation}

For a feasible point $\bar x$ of SR we formulate necessary optimality conditions with respect to the free variables from $I_1\left(\bar x\right)$:
\[
\frac{\partial f}{\partial x_j} \left( \bar x\right) = 0 \quad \mbox{for all } j\in I_1\left(\bar x\right).
\]
Recalling $f(x)=\frac{1}{2}\left\|Ax-b\right\|^2_2$, we get:
\[
\left(A^TA\bar x-A^Tb\right)_j=0 \quad \mbox{for all } j\in I_1\left(\bar x\right).
\]
Due to $\bar x_{I_0\left(\bar x\right)}=0$, it holds equivalently:
\begin{equation}
\label{eq:ts-cs}
A_{I_1\left(\bar x\right)}^TA_{I_1\left(\bar x\right)}\bar x_{I_1\left(\bar x\right)}-A_{I_1\left(\bar x\right)}^Tb=0.
\end{equation}
Note that $A_{I_1\left(\bar x\right)}$ stands for 
the submatrix of $A$ with the columns corresponding to the set $I_1\left(\bar x\right)$. Analogously, $x_{I_1\left(\bar x\right)}$ stands for the subvector of $x$ with the entries corresponding to the set $I_1\left(\bar x\right)$.

The previous derivation gives rise to the following definition. 

\begin{definition}[M-stationarity, \cite{burdakov:2016}]
\label{def:t-stat}
A feasible point $\bar x \in \R^{n,s}$ is called M-stationary for SR if 
\[
   A_{I_1\left(\bar x\right)}^TA_{I_1\left(\bar x\right)}\bar x_{I_1\left(\bar x\right)}-A_{I_1\left(\bar x\right)}^Tb=0.
\]
\end{definition}
Obviously, a local minimizer of SR is an M-stationary point.

Let us check the second-order sufficient optimality condition with respect to the free variables from $I_1\left(\bar x\right)$. We have:
\[
\displaystyle \left(\frac{\partial^2 f}{\partial x_j \partial x_k}\left(\bar x\right) \right)_{j,k \in {I_1\left(\bar x\right)}}  =  \displaystyle  A_{I_1\left(\bar x\right)}^TA_{I_1\left(\bar x\right)}.
\]
For the latter matrix to be positive definite, it is enough to assume that $A_{I_1\left(\bar x\right)}$ has full rank. Further, we examine the first-order behavior of $f$ on the sparse variables from $I_0\left(\bar x\right)$:
\[
  \frac{\partial f}{\partial x_i} \left( \bar x\right) = \left(A^TA\bar x-A^Tb\right)_i = \left(A_{I_0\left(\bar x\right)}^TA_{I_1\left(\bar x\right)}\bar x_{I_1\left(\bar x\right)}-A_{I_0\left(\bar x\right)}^Tb \right)_i \quad \mbox{for all } i\in I_0\left(\bar x\right).
\]
The following definition of nondegeneracy additionally requires the derivatives of the SR objective function with respect to the sparse variables be non-vanishing.

\begin{definition}[Nondegeneracy]
\label{def:nondeg}
An M-stationary point $\bar x \in \R^{n,s}$ of SR is called nondegenerate if the following conditions hold:
\begin{itemize}
\item[] ND1: if $\left\| \bar x\right\|_0 < s$ then all entries of the vector $A_{I_0\left(\bar x\right)}^TA_{I_1\left(\bar x\right)}\bar x_{I_1\left(\bar x\right)}-A_{I_0\left(\bar x\right)}^Tb$ are non-vanishing,
\item[] ND2: the matrix $A_{I_1\left(\bar x\right)}$ is of full rank, i.\,e. $\mbox{rank} \left(A_{I_1\left(\bar x\right)}\right)=\left\|\bar x\right\|_0$.
\end{itemize}
Otherwise, we call $\bar x$ degenerate.
\end{definition}

We point out that nondegeneracy is closely related to the property of $s$-regularity of the matrix $A$ introduced in \cite{beck:2013}.

\begin{definition}[$s$-regularity, \cite{beck:2013}]
A matrix $A\in \R^{m \times n}$ is called $s$-regular if for every index set
$I \subset \{1, \ldots ,n\}$ with $|I| = s$ it holds: $\mbox{rank}\left(A_I\right)=s$.
\end{definition}

\begin{lemma}[$s$-regularity and ND2]
\label{lem:sr-nd2}
If $A$ is $s$-regular, then ND2 is satisfied at all M-stationary points of SR. 
\end{lemma}

\proof
Let $\bar x$ be an M-stationary point of SR. Since we have $\left|I_1\left(\bar x\right)\right| \leq s$,  
there exists an index set
$I \subset \{1, \ldots ,n\}$ with $|I| = s$ and $I_1\left(\bar x \right) \subset I$. The $s$-regularity of $A$ implies that $\mbox{rank}\left(A_I\right)=s$. In particular, it follows that $\mbox{rank}\left(A_{I_1\left(\bar x \right)}\right)=\left\|\bar x\right\|_0$. \qed

\begin{lemma}[$s$-regularity and finiteness, \cite{beck:2013}]
\label{lem:sr-fin}
If $A$ is $s$-regular, then there are finitely many M-stationary points of SR.
\end{lemma}

\proof
If $\bar x$ is an M-stationary point of SR, then by using (\ref{eq:ts-cs}) we have:
\[
  \bar x_{I_1\left(\bar x\right)}=\left(A_{I_1\left(\bar x\right)}^TA_{I_1\left(\bar x\right)}\right)^{-1}A_{I_1\left(\bar x\right)}^Tb \quad \mbox{and} \quad  \bar x_{I_0\left(\bar x\right)}=0,
\]
where the matrix $A_{I_1\left(\bar x\right)}^TA_{I_1\left(\bar x\right)}$ is nonsingular due to the $s$-regularity of $A$. Since $\left|I_1\left(\bar x\right)\right| \leq s$, and the number of subsets of $\{1,2,\ldots,n\}$ with at most $s$ elements is finite, the result follows. \qed

Conditions ND1 and ND2 from Definition \ref{def:nondeg} allow to derive a relatively simple local representation of SR around a nondegenerate M-stationary point. In comparison to the corresponding result by \cite{laemmel:2019} for the sparsity  constrained nonlinear optimization (\ref{eq:scno}), the so-called quadratic index is vanishing here. This leads to the absence of negative squares in the representation (\ref{eq:normal}).

\begin{theorem}[Morse-Lemma for SR]
\label{thm:morse}
Suppose that $\bar x$ is a nondegenerate M-stationary point of SR. Then, there exist neighborhoods $U_{\bar x}$ and $V_0$ of $\bar x$ and $0$, respectively, and a local coordinate system $\Psi: U_{\bar x} \rightarrow V_0$ of $\R^n$ around $\bar x$ such that:

\begin{equation}
\label{eq:normal}
    f \circ \Psi^{-1}(y)=f\left(\bar x\right) + \sum\limits_{i \in I_0\left(\bar x\right)} y_i + \sum\limits_{j \in I_1\left(\bar x\right)} y_j^2,  
\end{equation}
where $y \in \R^{n-\NNorm{\bar x}{0},s-\NNorm{\bar x}{0}}\times \R^{\NNorm{\bar x}{0}}$.
\end{theorem}

\proof
Let $\bar x$ be a nondegenerate M-stationary point of SR.
By using the linear coordinate transformation $\Phi$ from (\ref{eq:diff}), we put $\bar f:= f\circ \Phi^{-1}$ on the set $\left(\R^{n-\NNorm{\bar x}{0},s-\NNorm{\bar x}{0}}\times \R^{\NNorm{\bar x}{0}} \right)\cap V_0$. 
As new coordinates we put $y=\left(y_{I_0\left(\bar x\right)},y_{I_1\left(\bar x\right)}\right)$. 
Then, it holds:
\[
 \begin{array}{rcl}
     \displaystyle \frac{\partial \bar f}{\partial y_{i}(0)} &=& \displaystyle \left(A_{I_0\left(\bar x\right)}^T A_{I_1\left(\bar x\right)}\bar x_{I_1\left(\bar x\right)}-A_{I_0\left(\bar x\right)}^Tb\right)_i \quad \mbox{for all } i \in I_0\left(\bar x\right),  \\ \\
     \displaystyle \frac{\partial \bar f}{\partial y_{j}(0)} &=&   \displaystyle \left( A_{I_1\left(\bar x\right)}^TA_{I_1\left(\bar x\right)}\bar x_{I_1\left(\bar x\right)}-A_{I_1\left(\bar x\right)}^Tb\right)_j \quad \mbox{for all } j \in I_1\left(\bar x\right), \\ \\
     \displaystyle \left(\frac{\partial^2 \bar f}{\partial y_j \partial y_k}(0) \right)_{j,k \in {I_1\left(\bar x\right)}}  &=&  \displaystyle  A_{I_1\left(\bar x\right)}^TA_{I_1\left(\bar x\right)}.
 \end{array}
\]
Due to ND1, M-stationarity of $\bar x$, and ND2, respectively, we have:
\begin{itemize}
    \item[(i)] if $\NNorm{\bar x}{0} < s$ then $\displaystyle \frac{\partial \bar f}{\partial y_i}(0) \ne 0$ for all $i \in {I_0\left(\bar x\right)}$,
    \item[(ii)] $\displaystyle \frac{\partial \bar f}{\partial y_j}(0) = 0$  for all $j\in {I_1\left(\bar x\right)}$,
\item[(iii)] the matrix $\displaystyle \left(\frac{\partial^2 \bar f}{\partial y_j \partial y_k}(0) \right)_{j,k \in {I_1\left(\bar x\right)}}$ is positive definite.
\end{itemize}
In what follows, we denote $\bar f$ by $f$ again. Under the following coordinate transformations the set $\R^{n-\NNorm{\bar x}{0},s-\NNorm{\bar x}{0}}\times \R^{\NNorm{\bar x}{0}}$ will be equivariantly transformed in itself. It holds:
\[
f\left(y\right)=\int_0^1 \frac{d}{dt} f\left(ty_{I_0\left(\bar x\right)},y_{I_1\left(\bar x\right)}\right)dt + f\left(0,y_{I_1\left(\bar x\right)}\right)=  \sum_{i \in I_0\left(\bar x\right)} y_i d_i(y)+f\left(0,y_{I_1\left(\bar x\right)}\right),
\]
whith linear functions $d_i$, $i \in I_0\left(\bar x\right)$.

Due to (ii)-(iii), we may apply the standard Morse lemma on the quadratic function $ f\left(0,y_{I_1\left(\bar x\right)}\right)$ without affecting the coordinates $y_{I_0\left(\bar x\right)}$, see e.\,g. \cite{Jongen:2000}. The corresponding coordinate transformation is linear.
Denoting the transformed functions again by $f$ and $d_i$, we obtain
\[
f(y) = f\left(\bar x\right)+\sum_{ i \in I_0\left(\bar x\right)}y_id_i(y) + \sum\limits_{j \in I_1\left(\bar x\right)} y_j^2.
\]

In case $\NNorm{\bar x}{0}=s$, we need to consider $f$ locally around the origin on the set 
\[
\R^{n-\NNorm{\bar x}{0},s-\NNorm{\bar x}{0}}\times \R^{\NNorm{\bar x}{0}} = \R^{n-\NNorm{\bar x}{0},0}\times \R^{\NNorm{\bar x}{0}} = \{0\}^{n-\NNorm{\bar x}{0}}\times \R^{\NNorm{\bar x}{0}}.
\]
Hence, $y_i=0$ for $i \in I_0\left(\bar x\right)$, and we immediately obtain the representation (\ref{eq:normal}).

In case $\NNorm{\bar x}{0}<s$, (i) provides that $\displaystyle d_i(0)=\frac{\partial f}{\partial y_i}(0) \not = 0$ for $i \in I_0\left(\bar x\right)$.
Hence, we may take
\[
 y_i d_i(y), i \in I_0\left(\bar x\right), \quad y_j,  j \in I_1\left(\bar x\right)
\]
as new local coordinates by a straightforward application of the inverse function theorem. Denoting the transformed function again by $f$, we obtain (\ref{eq:normal}).
Here, the coordinate transformation $\Psi$ is understood as the composition of all previous ones.
\qed

By means of Theorem \ref{thm:morse} the following important result follows. 



\begin{proposition}[Nondegenerate minimizers]
\label{prop:nondeg-min}
Let $\bar x$ be a nondegenerate M-stationary point for SR. Then, $\bar x$ is a local minimizer for SR if and only if the sparsity constraint is active, i.\,e. $\left\|\bar x\right\|_0=s$.
\end{proposition}

\proof
Let $\bar x$ be a nondegenerate M-stationary point of SR. The application of Morse Lemma from Theorem \ref{thm:morse} says that there exist neighborhoods $U_{\bar x}$ and $V_0$ of $\bar x$ and $0$, respectively, and a local $C^\infty$-coordinate system $\Psi: U_{\bar x} \rightarrow V_0$ of $\R^n$ around $\bar x$ such that:
\begin{equation}
\label{eq:help-nm0}
   f \circ \Psi^{-1}(y)=f\left(\bar x\right) + \sum\limits_{i \in I_0\left(\bar x\right)} y_i + \sum\limits_{j \in I_1\left(\bar x\right)} y_j^2,  
\end{equation}
where $y \in \R^{n-\NNorm{\bar x}{0},s-\NNorm{\bar x}{0}}\times \R^{\NNorm{\bar x}{0}}$. Therefore, $\bar x$ is a local minimizer for SR if and only if $0$ is a local minimizer of $f\circ \Psi^{-1}$ on the set $\left(\R^{n-\NNorm{\bar x}{0},s-\NNorm{\bar x}{0}}\times \R^{\NNorm{\bar x}{0}}\right) \cap V_0$. 
If 
we have $\NNorm{\bar x}{0}=s$, the formula in (\ref{eq:help-nm0}) reads as
\begin{equation}
\label{eq:help-nm1}
f\circ \Psi^{-1}(y)= f\left(\bar x\right) + \sum\limits_{j \in I_1\left(\bar x\right)} y_j^2,    
\end{equation}
where $y \in \{0\}^{n-s}\times \R^s$. Thus, $0$ is a local minimizer for (\ref{eq:help-nm1}). Vice versa, if $0$ is a local minimizer for (\ref{eq:help-nm0}), then obviously $\NNorm{\bar x}{0}=s$.
\qed

\section{Genericity}
\label{sec:gen}

Let us show that $s$-regularity is likely to be satisfied in the context of compressed sensing. This issue has  been already mentioned in \cite{beck:2013}.

\begin{lemma}[Genericity of $s$-regularity]
\label{lem:gen-sr}
Let $\mathcal{A}$ denote the subset of $s$-regular matrices $A$. Then, $\mathcal{A}$ is open and dense in $\R^{m \times n}$.
\end{lemma}
\proof
We consider the sets
\[
\Gamma_{I,r}=\left\{A \in \R^{m\times n} \; \left\vert \; \mbox{rank} \left(A_I\right)=r\right.\right\},
\]
where $I \subset\{1,\ldots,n\}$ with $|I|=s$, and $r=0,1,\ldots,s$. According to Example 7.3.23 from \cite{Jongen:2000}, $\Gamma_{I,r}$ is a submanifold of $\R^{m\times n}$ with codimension $(m-r)(s-r)$ -- recall that we have $r \leq s \leq m$ by Assumption \ref{ass:cs}.
In other words, $\Gamma_{I,r}$ is generically empty for $r=0,1,\ldots,s-1$, and $\Gamma_{I,s}$ is dense in $\R^{m\times n}$. Thus, $\mbox{rank}\left(A_I\right)=s$ holds for all $I$ in generic sense, which provides the assertion. \qed

Next, we show that ND1 and ND2 are fulfilled at all M-stationary points of SR for almost all 
data $(A,b) \in \R^{m \times n} \times \R^m$ with respect to Lebesgue measure, i.\,e. they are generically nondegenerate. 

\begin{theorem}[Genericity of nondegeneracy]
\label{thm:generic-cs}
Let $\mathcal{D}$ denote the subset of 
data $(A,b)$ for which each 
M-stationary point of SR is nondegenerate. Then, $\mathcal{D}$ is open and dense in $\R^{m \times n} \times \R^m$.
\end{theorem}

\proof
Due to Lemma \ref{lem:gen-sr}, the set $\mathcal{A}$ of $s$-regular matrices is open and dense in $\R^{m \times n}$. Then, Lemma \ref{lem:sr-nd2} implies that ND2 generically holds.
Now, we prove that ND1 is a generic condition for all M-stationary points $\bar x$ with $\NNorm{\bar x}{0}<s$. By setting $S=I_1\left(\bar x\right)$, we write (\ref{eq:ts-cs}) as  
\[
\bar x_S=\left(A_S^TA_S\right)^{-1}A_S^Tb=A_S^+b,
\]
where $A_S^+$ denotes the Moore-Penrose inverse of $A_S$. Hence, the vector under consideration in ND1 becomes:
\begin{equation}
\label{eq:h2}
    A_{S^c}^TA_S \bar x_S-A_{S^c}^Tb = - A_{S^c}^T\left(I - A_SA_S^+\right)b=- \left(\left(I - A_SA_S^+\right) A_{S^c}\right)^Tb,
\end{equation}
where we use the identity matrix $I \in \R^{m \times m}$, and the fact that $A_SA_S^+$ is symmetric:
\[
\left(A_SA_S^+\right)^T
     =\left(A_S^+\right)^TA_S^T 
     =\left(A_S^T\right)^+A_S^T
     =A_S\left(A_S^TA_S\right)^{-1}A_S^T
    = A_SA_S^+.
\]
Note that the entries of the vector in (\ref{eq:h2}) have to be shown generically non-vanishing, i.\,e.
\[
  \left(\left(I - A_SA_S^+\right) A_{\{i\}}\right)^Tb \not =0  \mbox{ for all } i \in S^c.
\]
%
%
For that, we define the sets
\[
\Omega_{S,i}=\left\{A\in\mathcal{A}\; \left\vert \; 
\left( I - A_SA_S^+\right)A_{\{i\}}=0 \right. \right\},
\]
where $S \subset \{1,\ldots,n\}$ with $|S| < s$, and $i \in S^c$.
Let us show that $\Omega_{S,i}$ is a submanifold. The condition
$\left( I - A_SA_S^+\right)A_{\{i\}}=0$
means that the vector $A_{\{i\}}$ lies in the nullspace of  $I - A_SA_S^+$, i.\,e. $A_{\{i\}} \in N\left( I - A_SA_S^+\right)$. 
Let us determine the dimension of $N\left( I - A_SA_S^+\right)$.
We start with the matrix $A_SA_S^+$.
It holds for the latter: 
\[
    \left(A_SA_S^+\right)A_SA_S^+ =
     A_S\left(A_S^TA_S\right)^{-1}A_S^TA_S\left(A_S^TA_S\right)^{-1}A_S^T = A_S\left(A_S^TA_S\right)^{-1}A_S^T=A_SA_S^+.
\]
Furthermore, the Sylvester's rank inequality provides
\[
\mbox{rank}\left(A_S A_S^+\right) \ge \mbox{rank}\left(A_S\right) + \mbox{rank}\left(A_S^+\right) - |S| = |S|,
\]
due to $\mbox{rank}\left(A_S^+\right)=\mbox{rank}\left(A_S\right)=|S|$ and the $s$-regularity of $A$.
Additionally, we have:
\[
\mbox{rank}\left(A_S A_S^+\right) \le \min\left\{\mbox{rank}\left(A_S \right),\mbox{rank}\left(A_S^+\right)\right\}=|S|.
\]
Altogether, $A_S A_S^+$ is an orthogonal projection
with $\mbox{rank}\left(A_S A_S^+\right)=|S|$.
Hence, $I-A_S A_S^+$ is also an orthogonal projection, and for the dimension of its nullspace we have: \[
\mbox{dim}\left(N\left( I - A_SA_S^+\right)\right)=\mbox{rank}\left(A_S A_S^+\right)=|S|.
\]
Since the dimension of $N\left( I - A_SA_S^+\right)$ remains 
constant under sufficiently small perturbations of $A_S$,  the condition $A_{\{i\}} \in N\left( I - A_SA_S^+\right)$ provides exactly $m-|S|$ stable equations. Hence, $\Omega_{S,i}$ is a submanifold of codimension $m-|S|$. Due to $m-|S| > m -s \geq 0$, we conclude that all $\Omega_{S,i}$ are generically empty.
In particular, the submanifold 
\[
\Omega_{S}=\left\{A\in\mathcal{A}\; \left\vert \; 
\left( I - A_SA_S^+\right)A_{\{i\}}\not = 0 \mbox{ for all } i\in S^c \right. \right\}
\]
is of dimension $mn$ and, hence, dense in $\R^{m\times n}$.

Finally, we define the sets
\[
\Upsilon_{S,i}=\left\{(A,b)\in \Omega_{S}\times \R^m\; \left\vert \; 
\left(\left( I - A_SA_S^+\right)A_{\{i\}}\right)^Tb = 0\right. \right\},
\]
where $S \subset \{1,\ldots,n\}$ with $|S| < s$, and $i \in S^c$ as above. Since $A \in \Omega_S$, the vector $\left( I - A_SA_S^+\right)A_{\{i\}}$ does not vanish. Hence, the equation $\left(\left( I - A_SA_S^+\right)A_{\{i\}}\right)^Tb = 0$ is nondegenerate. We conclude that $\Upsilon_{S,i}$ is a submanifold of codimension $1$, and can be therefore generically avoided. Overall, we have shown that the condition ND1 holds in generic sense.

The openness part follows due to the continuity of ND1 and ND2 with respect to sufficiently small perturbations of $A$.
\qed

We deduce the following important corollary on the structure of minimizers for SR.

\begin{corollary}[Sparsity constraint at minimizers]
\label{cor:min}
Generically, each minimizer $\bar x \in \R^{n,s}$ of SR is nondegenerate with the active sparsity constraint, i.\,e. $\left\|\bar x\right\|_0=s$.
\end{corollary}

\proof Note that every local minimizer of SR has to be M-stationary. Nondegenerate M-stationary points are generic by Theorem \ref{thm:generic-cs}. The rest follows by means of Proposition \ref{prop:nondeg-min}. \qed
%
%

\section{Stability}
\label{sec:stab}

Let us fix an arbitrary norm $\left\|(A,b)\right\|$ on the data space $(A,b) \in \R^{m \times n} \times \R^m$. For M-stationary points we define the notion of strong stability in the sense of \cite{kojima:1980}. 

\begin{definition}[Strong stability]
\label{def:s-stab}
An M-stationary point $\bar x$ of $\mbox{SR}(A,b)$ is called strongly stable if for some $r>0$ and each $\varepsilon \in (0,r]$ there exists $\delta >0$ such that whenever
\[
  \left(\widetilde A, \widetilde b \right) \in \R^{m \times n} \times \R^m \quad \mbox{and} \quad  \left\|\left(\widetilde A, \widetilde b \right) - (A,b)\right\| \leq \delta,
\]
the ball $B\left(\bar x, \varepsilon \right)$ contains an M-stationary point $\widetilde x$ of $\mbox{SR}\left(\widetilde A, \widetilde b \right)$ that is unique within the ball $B\left(\bar x, r \right)$.
\end{definition}

Let us illustrate a possible failure of strong stability of M-stationary points caused by their degeneracy.

\begin{example}[Instability]
\label{ex:cw-inst}
Let the following sensing matrix and measurement vector be given:
\[
   A= \left( \begin{array}{cc}
    1   &  0\\
    0   &  1
   \end{array}\right), \quad b = \left( \begin{array}{c}
    0\\0
   \end{array}\right).
\]
We consider the corresponding sparse  recovery problem with $s=1$:
\[
\mbox{SR}(A,b):\quad \min_{x_1,x_2}\,\, \frac{1}{2}x_1^2 + \frac{1}{2}x_2^2 \quad \mbox{s.\,t.} \quad 
   \left\|\left(x_1, x_2\right)\right\|_0 \leq 1.
\]
Obviously, $\bar x=(0,0)$ is the unique minimizer of $\mbox{SR}(A,b)$. Further, let us perturb the data by means of an arbitrarily small $\varepsilon > 0$ as follows:
\[
  \widetilde A = \left( \begin{array}{cc}
    1   &  0\\
    0   &  1
   \end{array}\right), \quad \widetilde b = \left( \begin{array}{c}
    \varepsilon \\ \varepsilon
   \end{array}\right).
\]
We obtain as perturbed sparse  recovery problem:
\[
\mbox{SR}\left(\widetilde A, \widetilde b\right):\quad \min_{x_1,x_2}\,\, \frac{1}{2}\left(x_1-\varepsilon\right)^2 + \frac{1}{2}\left(x_2-\varepsilon\right)^2 \quad \mbox{s.\,t.} \quad 
   \left\|\left(x_1, x_2\right)\right\|_0 \leq 1.
\]
It is easy to see that $\mbox{SR}\left(\widetilde A, \widetilde b\right)$ has now two solutions $\widetilde x^{1a}=(\varepsilon,0)$ and $\widetilde x^{1b}=(0,\varepsilon)$. 
Here, we observe a bifurcation of the minimum $\bar x$ of the original problem $\mbox{SR}(A,b)$ into two minima $\widetilde x^{1a}$ and $\widetilde x^{1b}$ of the perturbed problem $\mbox{SR}\left(\widetilde A, \widetilde b\right)$. Let us explain this bifurcation in terms of M-stationarity. The bifurcation is caused by the degeneracy of $\bar x$ viewed as an M-stationary point of $\mbox{SR}(A,b)$. Note that ND1 is violated at the M-stationary point $\bar x$ of $\mbox{SR}(A,b)$. More interestingly, there is another M-stationary point $\widetilde x^2 =(0,0)$ of the perturbed problem. In fact, due to $\left\| \widetilde x^2\right\|_0=0$ and the validity of ND1, $\widetilde x^2$ is a nondegenerate M-stationary point of $\mbox{SR}\left(\widetilde A, \widetilde b\right)$. For the latter we have 
\[
\left\|\widetilde x^2\right\|_0=s-1,
\]
meaning that $\widetilde x^2$ is a saddle point which connects two nondegenerate minimizers $\widetilde x^{1a}$ and $\widetilde x^{1b}$ of $\mbox{SR}\left(\widetilde A, \widetilde b\right)$. Overall, we conclude that the degenerate minimum $\bar x$ of the original problem $\mbox{SR}(A,b)$ is not strongly stable. Moreover, it bifurcates into two nondegenerate minima $\widetilde x^{1a}$ and $\widetilde x^{1b}$, as well as leads to one nondegenerate saddle point $\widetilde x^2$ of the perturbed problem $\mbox{SR}\left(\widetilde A, \widetilde b\right)$. It is not hard to see that 
every of the M-stationary points $\widetilde x^{1a}$, $\widetilde x^{1b}$, and 
$\widetilde x^2$ are strongly stable for $\mbox{SR}\left(\widetilde A, \widetilde b\right)$. \qed
\end{example}

It turns out that Example \ref{ex:cw-inst} is typical in the context of sparse  recovery. Namely, strong stability and nondegeneracy are equivalent properties of M-stationary points.

\begin{theorem}[Characterization of strong stability]
\label{thm:stab-nondeg}
An M-stationary point $\bar x$ of $\mbox{SR}(A,b)$ is strongly stable if and only if it is nondegenerate.
\end{theorem}
\proof We start with the {\it necessity part}. Let $\bar x \in \R^{n,s}$ be a nondegenerate M-stationary point of $\mbox{SR}(A,b)$. For any $\left(\widetilde A, \widetilde b \right)$ chosen sufficiently close to $(A,b)$, we show that there exists a unique M-stationary point $\widetilde x \in \R^{n,s}$ of SR$\left(\widetilde A, \widetilde b \right)$ in a neighborhood of $\bar x$. First, for all $\widetilde x \in \R^{n,s}$ being sufficiently close to $\bar x$ we have by continuity arguments that
\[
   I_1\left(\widetilde x \right) \supset I_1\left(\bar x\right).
\]
We claim that if $\widetilde x$ is additionally an M-stationary point of SR$\left(\widetilde A, \widetilde b \right)$ then actually the equality holds above, i.\,e.
\[
  I_1\left(\widetilde x \right) = I_1\left(\bar x\right).
\]
To see this we consider the following cases:
\begin{itemize}
    \item[(i)] The sparsity constraint is active, i.\,e. $\left\|\bar x\right\|_0=s$. Then, due to $\left\|\widetilde x\right\|_0\leq s$, we have trivially $I_1\left(\widetilde x \right) = I_1\left(\bar x\right)$.
    \item[(ii)] The sparsity constraint is not active, i.\,e. $\left\|\bar x\right\|_0<s$. We assume in contrary that there exists $\widetilde i \in I_1\left(\widetilde x \right) \backslash I_1\left(\bar x\right)$. By having $\widetilde i \in I_0\left( \bar x \right)$ and recalling ND1 for $\bar x$, we obtain
    \[
       A_{\widetilde i}^TA_{I_1\left(\bar x\right)}\bar x_{I_1\left(\bar x\right)}-A_{\widetilde i}^Tb  \not = 0.
    \]
By continuity, it follows that
    \[
       \widetilde A_{\widetilde i}^T \widetilde A_{I_1(\widetilde x)}\widetilde x_{I_1(\widetilde x)}-\widetilde A_{\widetilde i}^T \widetilde b  \not = 0
    \]
for any $\left(\widetilde A, \widetilde b \right)$ from a sufficiently small neighborhood of $(A,b)$. Due to $\widetilde i \in I_1\left( \widetilde x \right)$,
the latter contradicts the fact that $\widetilde x$ is M-stationary for SR$\left(\widetilde A, \widetilde b \right)$. 
\end{itemize}
Further, for $\widetilde x \in \R^{n,s}$ to be M-stationary for SR$\left(\widetilde A, \widetilde b \right)$ the following holds, see Definition \ref{def:t-stat}:
  \[
   \widetilde A_{I_1\left(\widetilde x\right)}^T \widetilde A_{I_1\left(\widetilde x\right)}\widetilde  x_{I_1\left(\widetilde x\right)}-\widetilde A_{I_1\left(\widetilde x\right)}^T \widetilde b=0.
\]
Since $I_1\left(\widetilde x \right) = I_1\left(\bar x\right)$, we get locally:
\[
 \widetilde A_{I_1(\bar  x)}^T \widetilde A_{I_1(\bar  x)}\widetilde  x_{I_1(\bar  x)}-\widetilde A_{I_1(\bar  x)}^T \widetilde b=0.
\]
Moreover, by continuity and ND2 for $\bar x$, the matrix $\widetilde A_{I_1(\bar  x)}$
is of full rank, i.\,e. $\mbox{rank}\left(\widetilde A_{I_1\left(\bar x\right)}\right)=\left\|\bar x\right\|_0$. Consequently, the unique M-stationary point of SR$\left(\widetilde A, \widetilde b \right)$ in a neighborhood of $\bar x$ is
\[
\widetilde  x_{I_1(\bar  x)} = \left(\widetilde A_{I_1(\bar  x)}^T \widetilde A_{I_1(\bar  x)} \right)^{-1} \widetilde A_{I_1(\bar  x)}^T \widetilde b \quad \mbox{and} \quad
 \widetilde x_{I_0(\bar  x)} =0,
\]
which depends continuously on $\left(\widetilde A, \widetilde b \right)$.

We proceed with the {\it sufficiency part}. Let $\bar x \in \R^{n,s}$ be a strongly stable M-stationary point of $\mbox{SR}(A,b)$. We show by contradiction that $\bar x$ is also nondegenerate.

First, let us assume that ND2 is violated for $\bar x$, hence, the matrix $A_{I_1\left(\bar x\right)}$ is not of full rank, i.\,e. $\mbox{rank} \left(A_{I_1\left(\bar x\right)}\right)<\left\|\bar x\right\|_0$. We consider the following system of linear equations with respect to $x$:
\begin{equation}
  \label{eq:lsys} 
  A_{I_1\left(\bar x\right)}^TA_{I_1\left(\bar x\right)} x_{I_1\left(\bar x\right)}-A_{I_1\left(\bar x\right)}^T b=0 \quad \mbox{and} \quad x_{I_0\left(\bar x\right)}=0.
\end{equation}
Note that $\bar x$ solves (\ref{eq:lsys}) as an M-stationary point of $\mbox{SR}(A,b)$.   Since the matrix $A_{I_1\left(\bar x\right)}^TA_{I_1\left(\bar x\right)}$ is singular, the solution set of (\ref{eq:lsys}) is a linear subspace of dimension $\left\|\bar x\right\|_0 -\mbox{rank} \left(A_{I_1\left(\bar x\right)}\right) >0$.
Any solution $x$ of (\ref{eq:lsys}) is feasible for $\mbox{SR}(A,b)$, since $\left\|x\right\|_0 \leq \left\|\bar x\right\|_0\leq s$.
Moreover, within a sufficiently small neighborhood of $\bar x$ it holds
$x_i\not = 0$ for all $i \in I_1\left(\bar x\right)$, hence, $I_1\left(x\right)=I_1\left(\bar x\right)$.
Altogether, the solutions of (\ref{eq:lsys}) around $\bar x$ are M-stationary for $\mbox{SR}(A,b)$. Thus, $\bar x$ is not isolated as an M-stationary point of $\mbox{SR}(A,b)$ and, therefore, cannot be strongly stable, a contradiction.

Second, we assume that ND1 is violated for $\bar x$, but ND2 is fulfilled. Then, we have $\left\| \bar x\right\|_0 < s$ and there exists $\bar i \in I_0\left(\bar x\right)$ with
\begin{equation}
  \label{eq:lshelp} 
       A_{\bar i}^TA_{I_1\left(\bar x\right)}\bar x_{I_1\left(\bar x\right)}-A_{\bar i}^Tb  = 0.
\end{equation}
As an auxiliary claim, we show that the matrix $A_{I_1\left(\bar x\right) \cup \left\{\bar i\right\}}$ is of full rank, i.\,e. 
\begin{equation}
\label{eq:auxil}
\mbox{rank} \left(A_{I_1\left(\bar x\right) \cup \left\{\bar i\right\}}\right)= \left\|\bar x\right\|_0+1. 
\end{equation}
Let us assume for a moment that (\ref{eq:auxil}) is not fulfilled. We come to a contradiction by considering the following system of linear equations with respect to $x$:
\begin{equation}
  \label{eq:lsys3} 
  A_{I_1\left(\bar x\right) \cup \left\{\bar i\right\}}^T A_{I_1\left(\bar x\right) \cup \left\{\bar i\right\}} x_{I_1\left(\bar x\right) \cup \left\{\bar i\right\}}- A_{I_1\left(\bar x\right) \cup \left\{\bar i\right\}}^T b=0 \quad \mbox{and} \quad x_{I_0\left(\bar x\right) \backslash \left\{\bar i\right\}}=0.
\end{equation}
Note that $\bar x$ solves (\ref{eq:lsys3}) as an M-stationary point of $\mbox{SR}(A,b)$ and due to (\ref{eq:lshelp}). Since we suppose that $\mbox{rank} \left(A_{I_1\left(\bar x\right) \cup \left\{\bar i\right\}}\right) < \left\|\bar x\right\|_0+1$, the matrix $A_{I_1\left(\bar x\right) \cup \left\{\bar i\right\}}^T A_{I_1\left(\bar x\right) \cup \left\{\bar i\right\}}$ is singular, and the solution set of (\ref{eq:lsys3}) is a linear subspace of dimension 
\[
\left\|\bar x\right\|_0+1 -\mbox{rank} \left(A_{I_1\left(\bar x\right) \cup \left\{\bar i\right\}}\right) >0.
\]
Any solution $x$ of (\ref{eq:lsys3}) is feasible for $\mbox{SR}(A,b)$, since $\left\|x\right\|_0 \leq \left\|\bar x\right\|_0 +1 \leq s$.
Moreover, within a sufficiently small neighborhood of $\bar x$ it holds
$x_i\not = 0$ for all $i \in I_1\left(\bar x\right)$. 
Let us assume for a moment that for all solutions $x$ of (\ref{eq:lsys3}) within an arbitrarily small  neighborhood of $\bar x$, it holds $x_{\bar i} = 0$. Then, $x$ solves (\ref{eq:lsys}). Due to ND2, the matrix $A_{I_1\left(\bar x\right)}^TA_{I_1\left(\bar x\right)}$ is nonsingular, which implies that $x=\bar x$. This would mean that (\ref{eq:lsys3}) is uniquely solvable in a neighborhood of $\bar x$, a contradiction to the singularity of $A_{I_1\left(\bar x\right) \cup \left\{\bar i\right\}}^T A_{I_1\left(\bar x\right) \cup \left\{\bar i\right\}}$. 
Altogether, in any sufficiently small neighborhood of $\bar x$ there exist solutions $x$ of (\ref{eq:lsys3}) such that $I_1\left(x\right)=I_1\left(\bar x\right) \cup \left\{\bar i\right\}$. Hence, those points $x$ are M-stationary for $\mbox{SR}(A,b)$. Thus, $\bar x$ is not isolated as an M-stationary point of $\mbox{SR}(A,b)$ and, therefore, cannot be strongly stable. By contradiction, we just conclude that (\ref{eq:auxil}) is fulfilled. 

After this preliminary considerations, we construct a perturbation $\left(\widetilde A, \widetilde b\right)$ arbitrarily close to $(A,b)$ such that $\mbox{SR}\left(\widetilde A, \widetilde b\right)$ has at least two M-stationary points $\widetilde x^1 \not = \widetilde x^2$ in a proximity to $\bar x$. We may assume that for all $\widetilde b$ close to $b$ it holds:
\[
\widetilde b \not = A_{I_1\left(\bar x\right)} \bar x_{I_1\left(\bar x\right)}.
\]
In fact, if $b = A_{I_1\left(\bar x\right)} \bar x_{I_1\left(\bar x\right)}$ then any $\widetilde b \not = b$ suffices. 
Hence, there exists a normalized vector $c \in \R^m$ such that for all $\widetilde b$ close, but not necessarily equal to $b$ it holds:
\[
   c^T\left( A_{I_1\left(\bar x\right)} \bar x_{I_1\left(\bar x\right)} - \widetilde b \right) \not = 0.
\]
We define the following family of perturbations $\widetilde A \in \R^{m,n}$ depending on the parameter $t\in \R^2$:
\[
  \widetilde A_{\{\bar i\}}=  A_{\{\bar i\}} + t c, \quad
  \widetilde A_{\{\bar i\}^c} = A_{\{\bar i\}^c}.
\]
Note that $\widetilde A$ differs from $A$ only with respect to the $\bar i$-th column. Moreover, for $t=0$ both matrices $\widetilde A$ and $A$ coincide.


{\it 1) Construction of $\widetilde x^1$.} We consider the following system of linear equations with respect to $x$:
\begin{equation}
  \label{eq:lsys1} 
  \widetilde A_{I_1\left(\bar x\right)}^T \widetilde A_{I_1\left(\bar x\right)} x_{I_1\left(\bar x\right)}-\widetilde A_{I_1\left(\bar x\right)}^T \widetilde b=0 \quad \mbox{and} \quad x_{I_0\left(\bar x\right)}=0.
\end{equation}
Since $\bar i \not \in I_1\left(\bar x\right)$, the system (\ref{eq:lsys1}) is equivalent to
\[
A_{I_1\left(\bar x\right)}^T A_{I_1\left(\bar x\right)} x_{I_1\left(\bar x\right)}- A_{I_1\left(\bar x\right)}^T \widetilde b=0 \quad \mbox{and} \quad x_{I_0\left(\bar x\right)}=0.
\]
Due to ND2, the unique solution of (\ref{eq:lsys1}) is then
\[
\widetilde x^1_{I_1\left(\bar x\right)}= \left( A_{I_1\left(\bar x\right)}^T  A_{I_1\left(\bar x\right)}\right)^{-1}  A_{I_1\left(\bar x\right)}^T \widetilde b \quad \mbox{and} \quad \widetilde x^1_{I_0\left(\bar x\right)}=0.
\]
Note that $\widetilde x^1$ is independent of $t$. We see that $\widetilde x^1$ is feasible for $\mbox{SR}\left(\widetilde A, \widetilde b\right)$, since $\left\|\widetilde x^1\right\|_0 \leq \left\|\bar x\right\|_0<s$. Since $\widetilde x^1$ depends continuously on $\widetilde b$, the point
$\widetilde x^1$ falls into a sufficiently small neighborhood of $\bar x$.
Hence, it holds $\widetilde x^1_i \not = 0$ for all $i \in I_1\left(\bar x\right)$ or, equivalently, $I_1\left(\widetilde x^1\right)=I_1\left(\bar x\right)$. Altogether, $\widetilde x^1$ is M-stationary for $\mbox{SR}\left(\widetilde A, \widetilde b\right)$.

{\it 2) Construction of $\widetilde x^2$.} We consider the following system of linear equations with respect to $x$:
\begin{equation}
  \label{eq:lsys2} 
  \widetilde A_{I_1\left(\bar x\right) \cup \left\{\bar i\right\}}^T \widetilde A_{I_1\left(\bar x\right) \cup \left\{\bar i\right\}} x_{I_1\left(\bar x\right) \cup \left\{\bar i\right\}}-\widetilde A_{I_1\left(\bar x\right) \cup \left\{\bar i\right\}}^T \widetilde b=0 \quad \mbox{and} \quad x_{I_0\left(\bar x\right) \backslash \left\{\bar i\right\}}=0.
\end{equation}
As we showed in (\ref{eq:auxil}), the matrix $A_{I_1\left(\bar x\right) \cup \left\{\bar i\right\}}$ is of full rank. Due to continuity reasons, the matrix $\widetilde A_{I_1\left(\bar x\right) \cup \left\{\bar i\right\}}$ is also of full rank at least for $t$ sufficiently close to zero. The unique solution of (\ref{eq:lsys2}) is then
\[
  \widetilde x^2_{I_1\left(\bar x\right) \cup \left\{\bar i\right\}} = \left( \widetilde A_{I_1\left(\bar x\right) \cup \left\{\bar i\right\}}^T \widetilde A_{I_1\left(\bar x\right) \cup \left\{\bar i\right\}}\right)^{-1} \widetilde A_{I_1\left(\bar x\right) \cup \left\{\bar i\right\}}^T \widetilde b \quad \mbox{and} \quad   \widetilde x^2_{I_0\left(\bar x\right) \backslash \left\{\bar i\right\}}=0.
\]
We see that $\widetilde x^2$ is feasible for $\mbox{SR}\left(\widetilde A, \widetilde b\right)$, since $\left\|\widetilde x^2\right\|_0 \leq \left\|\bar x\right\|_0 +1 \leq s$. 
Now, we use the fact that $\bar x$ is the unique solution of 
(\ref{eq:lsys3}), i.\,e.
\[
\bar x_{I_1\left(\bar x\right) \cup \left\{\bar i\right\}} = \left(  A_{I_1\left(\bar x\right) \cup \left\{\bar i\right\}}^T A_{I_1\left(\bar x\right) \cup \left\{\bar i\right\}}\right)^{-1} A_{I_1\left(\bar x\right) \cup \left\{\bar i\right\}}^T b \quad \mbox{and} \quad  \bar x_{I_0\left(\bar x\right) \backslash \left\{\bar i\right\}}=0.
\]
As consequence, $\widetilde x^2$ falls into an arbitrarily small neighborhood of $\bar x$ as soon as $\left(\widetilde A,\widetilde b \right)$ is sufficiently close to $(A,b)$. This implies that $\widetilde x^2_i \not = 0$ for all $i \in I_1\left(\bar x\right)$. 
In order to show that $\widetilde x^2_{\bar i} \not =0$, we compute the derivative of $\widetilde x^2_{I_1\left(\bar x\right) \cup \left\{\bar i\right\}}$ at $t=0$ by using the implicit function theorem for the system of equations (\ref{eq:lsys2}).
For its left-hand side we set
\[
  G\left(x_{I_1\left(\bar x\right) \cup \left\{\bar i\right\}},t\right) = \widetilde A_{I_1\left(\bar x\right) \cup \left\{\bar i\right\}}^T \widetilde A_{I_1\left(\bar x\right) \cup \left\{\bar i\right\}} x_{I_1\left(\bar x\right) \cup \left\{\bar i\right\}}-\widetilde A_{I_1\left(\bar x\right) \cup \left\{\bar i\right\}}^T \widetilde b.   
\]
It is straight-forward to see that
\[
   D_{x_{I_1\left(\bar x\right) \cup \left\{\bar i\right\}}}G\left(\bar x_{I_1\left(\bar x\right) \cup \left\{\bar i\right\}},0\right) = \widetilde A_{I_1\left(\bar x\right) \cup \left\{\bar i\right\}}^T \widetilde A_{I_1\left(\bar x\right) \cup \left\{\bar i\right\}}
\]
and
\[
  D_{t}G\left(\bar x_{I_1\left(\bar x\right) \cup \left\{\bar i\right\}},t\right) = \left( \begin{array}{c}
        0 \\ \vdots \\ 0 \\ c^T\left( A_{I_1\left(\bar x\right)} \bar x_{I_1\left(\bar x\right)} - \widetilde b \right) \\ 0 \\ \vdots \\ 0
  \end{array}{}\right),
\]
where $c^T\left( A_{I_1\left(\bar x\right)} \bar x_{I_1\left(\bar x\right)} - \widetilde b \right)$ is the $\bar i$-th component of $D_{t}G\left(\bar x_{I_1\left(\bar x\right) \cup \left\{\bar i\right\}},t\right)$.
It follows that
\begin{equation}
\label{eq:lshelp3}
\begin{array}{lcl}
     \displaystyle \frac{\mbox{d} \widetilde x^2_{I_1\left(\bar x\right) \cup \left\{\bar i\right\}}}{\mbox{d}t} &=& \displaystyle -  \left( D_{x_{I_1\left(\bar x\right) \cup \left\{\bar i\right\}}}G\left(\bar x_{I_1\left(\bar x\right) \cup \left\{\bar i\right\}},0\right) \right)^{-1}  D_{t}G\left(\bar x_{I_1\left(\bar x\right) \cup \left\{\bar i\right\}},t\right) \\ \\
   &=& - \displaystyle \left(\widetilde A_{I_1\left(\bar x\right) \cup \left\{\bar i\right\}}^T \widetilde A_{I_1\left(\bar x\right) \cup \left\{\bar i\right\}}\right)^{-1}\left( \begin{array}{c}
        0 \\ \vdots \\ 0 \\ c^T\left( A_{I_1\left(\bar x\right)} \bar x_{I_1\left(\bar x\right)} - \widetilde b \right) \\ 0 \\ \vdots \\ 0
  \end{array}{}\right).
\end{array}
\end{equation}
Let us assume for a moment that $\widetilde x^2_{\bar i} = 0$ for all $t$ within a neighborhood of zero. Then, every $\widetilde x^2$ solves also (\ref{eq:lsys1}), hence, $\widetilde x^2 = \widetilde x^1$.
In this case $\widetilde x^2_{I_1\left(\bar x\right) \cup \left\{\bar i\right\}}$ is constant and its derivative with respect to $t$ vanishes around zero. Substituting into (\ref{eq:lshelp3}), we obtain
\[
   c^T\left( A_{I_1\left(\bar x\right)} \bar x_{I_1\left(\bar x\right)} - \widetilde b \right) = 0,
\]
a contradiction to the choice of $\widetilde b$. 
We have just shown that in any sufficiently small neighborhood of $\bar x$ there exist solutions $\widetilde x^2$ of (\ref{eq:lsys2}) such that $I_1\left(\widetilde x^2\right)=I_1\left(\bar x\right) \cup \left\{\bar i\right\}$. Altogether, $\widetilde x^2$ is M-stationary for $\mbox{SR}\left(\widetilde A, \widetilde b\right)$.
We have also shown that $\widetilde x^1 \not = \widetilde x^2$, since $\widetilde x^1_{\bar i} =0$ and $\widetilde x^2_{\bar i} \not =0$. By this, the M-stationary point $\bar x$ of $\mbox{SR}(A,b)$ is not strongly stable. \qed

We point out that the equivalence of strong stability and nondegeneracy of stationary point is by far not usual in nonsmooth optimization. Exemplarily, let us compare the relation of both notions in the context of mathematical programs with complementarity constraints. 

\begin{example}[Complementarity constraints]
Let the following sensing matrix and measurement vector be given:
\[
   A= \left( \begin{array}{cc}
    1   &  0\\
    0   &  1
   \end{array}\right), \quad b = \left( \begin{array}{c}
    -1\\0
   \end{array}\right).
\]
We consider the corresponding sparse  recovery problem with $s=1$:
\[
\mbox{SR}(A,b):\quad \min_{x_1,x_2}\,\, \frac{1}{2}(x_1+1)^2 + \frac{1}{2}x_2^2 \quad \mbox{s.\,t.} \quad 
   \left\|\left(x_1, x_2\right)\right\|_0 \leq 1.
\]
Obviously, $(-1,0)$ is the nondegenerate minimizer of $\mbox{SR}(A,b)$. However, there exists another M-stationary point of $\mbox{SR}(A,b)$, namely, $\bar x=(0,0)$. Due to the violation of ND1, $\bar x$ is degenerate and, hence, cannot be strongly stable for $\mbox{SR}(A,b)$ in view of Theorem \ref{thm:stab-nondeg}. Now, we consider the following mathematical program with complementarity constraints:
\begin{equation}
\label{eq:cc}
\mbox{MPCC}(A,b):\quad \min_{x_1,x_2}\,\, f\left(x_1,x_2\right)=\frac{1}{2}(x_1+1)^2 + \frac{1}{2}x_2^2 \quad \mbox{s.\,t.} \quad 
   x_1 \cdot x_2=0, \quad x_1, x_2 \geq 0.
\end{equation}
The objective function $f$ of MPCC$(A,b)$ is the same as of $\mbox{SR}(A,b)$, but the sparsity constraint $\left\|\left(x_1, x_2\right)\right\|_0 \leq 1$ is substituted by the complemenarity constraint $x_1 \cdot x_2=0, x_1, x_2 \geq 0$. We see that $\bar x=(0,0)$ is the unique minimizer of MPCC$(A,b)$. In particular, $\bar x$ is a so-called C-stationary point of MPCC$(A,b)$, see e.\,g. \cite{jongen:2012} for details. In fact, at $\bar x$ the derivatives of the objective function with respect to biactive variables are of the same sign:
\[
  \frac{\partial f}{\partial x_1}\left(\bar x\right) =1, \quad 
  \frac{\partial f}{\partial x_2}\left(\bar x\right)=0.
\]
However, $\bar x$ is degenerate, since one of the above derivatives vanishes. Nevertheless, $\bar x$ is a strongly stable C-stationary point of MPCC$(A,b)$. This is due to Corollary 3.1 from \cite{jongen:2012}, where the following sufficient condition for the strong stability of C-stationary points is given:
\[
  \frac{\partial f}{\partial x_2}\left(\bar x\right)=0, \quad \frac{\partial^2 f}{\partial x^2_2}\left(\bar x\right) \cdot 
   \frac{\partial f}{\partial x_1}\left(\bar x\right) >0. 
\]
Here, the latter is fulfilled due to
\[
   \frac{\partial^2 f}{\partial x^2_2}\left(\bar x\right) = 1, \quad
   \frac{\partial f}{\partial x_1}\left(\bar x\right) =1.
\]
We conclude that for mathematical programs with complementarity constraints the equivalence of strong stability and nondegeneracy is not valid as it is the case of sparse  recovery. \qed
\end{example}

\section{Global aspects}
\label{sec:glob}

Let us study the topological properties of SR lower level sets
\[
M^{a}=\left\{ x \in\R^{n,s}\, \left\vert \, f(x)\le a \right. \right\},
\]
where $a \in \R$ is varying. For that, we define intermediate sets for $a<b$:
\[
M^{b}_{a}=\left\{ x \in\R^{n,s}\, \left\vert \, a \leq f(x) \leq b \right. \right\}.
\]
For the topological concepts used below we
refer to \cite{spanier:1966}.

In what follows, we mention several consequences of the $s$-regularity of $A$ for the topological properties of the lower level sets. 

\begin{lemma}[Lower level sets]
\label{lem:help}
Let $A$ be an $s$-regular matrix. Then, all lower level sets $M^a$ are bounded. Moreover, for all sufficiently large $a\in \R$ they are also connected.  
\end{lemma}

\proof
First, we show that the lower level sets 
\[
M^a = \left\{x \in \R^{n,s} \, \left | \, \left\|Ax-b\right\|^2_2 \leq a \right.\right\}
\]
are bounded for any $a \in \R$. We write for the SR feasible set:
\[
   \R^{n,s} = \bigcup_{\scriptsize \begin{array}{c}
        S \subset \{1, \ldots,n\} \\
        |S| \leq s 
   \end{array}} X_{S},
\]   
where 
\[
  X_{S} = \left\{x\in \R^n \,\left|\, x_{S^c} =0 \right.\right\}.
\]
Hence, 
\begin{equation}
    \label{eq:repr}  
        M^a = \bigcup_{\scriptsize \begin{array}{c}
        S \subset \{1, \ldots,n\} \\
        |S| \leq s 
   \end{array}} M_S^a,
\end{equation}
where
\[
 M_S^a=\left\{x \in X_S \, \left | \, \left\|Ax-b\right\|^2_2 \leq a \right.\right\}.
\]
It holds for $x \in X_S$:
\[
  \left\|Ax-b\right\|^2_2 = \left\|A_S x_S-b\right\|^2_2.
\]
Since $A$ is $s$-regular and $|S| \leq s$, we have $\mbox{rank}\left(A_S\right)=|S|$. Hence, the sets
\[
   M_S^a= \left\{x \in X_S \, \left | \, \left\|A_Sx_S-b\right\|^2_2 \leq a \right.\right\}
\]
are $|S|$-dimensional ellipsoids and, thus, bounded. Therefore, the set $M^a$ is bounded as well, namely, as a finite union of bounded sets $M_S^a$, cf. the representation (\ref{eq:repr}).

Further, it is possible to increase $a \in \R$ in order to guarantee that $M^a$ is also connected. To have this, let us assume that $a \geq  \|b\|_2^2$. Hence, we have $0 \in M_S^a$ for all $S \subset \{1, \ldots,n\}$ with $|S| \leq s$.
Moreover, the sets $M_S^a$ from the representation (\ref{eq:repr}) are connected as $|S|$-dimensional ellipsoids. As a consequence, the set $M^a$ is connected as a finite union of connected sets $M_S^a$, all of them having nonempty intersection, i.\,e.
\[
  0 \in \bigcap_{\scriptsize \begin{array}{c}
        S \subset \{1, \ldots,n\} \\
        |S| \leq s 
   \end{array}} M_S^a.
\]
This concludes the proof. \qed

We show that the lower level sets do not undergo topological changes when passing a non-M-stationary level.

\begin{theorem}[Deformation for SR]
\label{thm:def}
Let $A$ be an $s$-regular matrix and $M^b_a$ contain no M-stationary points for SR. Then, $M^a$ is homeomorphic to $M^b$.
\end{theorem}

\proof
We apply Proposition 3.2 from Part I in \cite{goresky:1988}. The latter provides the deformation for general Whitney stratified sets with respect to critical points of proper maps. Note that the SR feasible set admits a Whitney stratification:
\[
   \R^{n,s} = \bigcup_{\scriptsize \begin{array}{c}
        I \subset \{1, \ldots,n\} \\
        |I| \leq s 
   \end{array}} \bigcup_{J \subset I} Z_{I,J},
\]   
where 
\[
  Z_{I,J} = \left\{x\in \R^n \,\left|\, x_{I^c} =0, x_J > 0, x_{I\backslash J} <0 \right.\right\}.
\]
The notion of criticality used in \cite{goresky:1988} can be stated for SR as follows. A point $\bar x \in \R^{n,s}$ is called critical for $f$ on $\R^{n,s}$ if it holds:
\[
Df\left(\bar x\right)_{|T_{\bar x} Z} =0,
\]
where $Z$ is the stratum of $\R^{n,s}$ which contains $\bar x$, and $T_{\bar x} Z$ is the tangent space of $Z$ at $\bar x$. By identifying $I =I_1\left(\bar x\right)$ and, hence, $I^c =I_0\left(\bar x\right)$, we see that the concepts of criticality and M-stationarity coincide. It remains to note that, due to Lemma \ref{lem:help}, the restriction of $f$ on  $\R^{n,s}$ is proper, i.\,e. $f^{-1}(K)\cap \R^{n,s}$ is compact for any compact set $K \subset \R$. \qed

Now, we turn our attention to the topological changes of lower level sets when passing an M-stationary level. Traditionally, they are described by means of the so-called cell-attachment.
We first consider a special case of cell-attachment dealt with already in \cite{laemmel:2019}. For that, let $N^\epsilon$ denote the lower level set of a special linear function on $\R^{p,q}$, i.\,e.
\[
 N^\epsilon = \left\{ x \in \R^{p,q} \,\left|\, \sum_{i=1}^{p} x_i \leq \epsilon \right. \right\},
\]
where $\epsilon \in \R$, and the integers $q < p$ are nonnegative. 

\begin{lemma}[Normal Morse data, \cite{laemmel:2019}]
\label{lem:cat}
For any $\epsilon > 0$ the set $N^\epsilon$ is homotopy-equivalent to $N^{-\epsilon}$ with $\binom{p-1}{q}$ cells of dimension $q$ attached. The latter cells are the $q$-dimensional simplices from the collection 
\[
 \left\{ \left. \mbox{conv} \left(e_j, j \in J\right)  \,\right|\, J \subset \{1, \ldots,p\}, 1 \in J, |J| = q+1 \right\}.
\]
\end{lemma}

The general case of cell-attachment can be shown by using Lemma \ref{lem:cat}. 

\begin{theorem}[Cell-Attachment for SR]
\label{thm:cell-a}
Let $A$ be an $s$-regular matrix and
$M^b_a$ contain exactly one M-stationary point $\bar x$ for SR.
If $a<f\left(\bar x \right) <b$,
then $M^b$ is homotopy-equivalent to $M^a$ with $\binom{n-\left\|\bar x\right\|_0-1}{s-\left\|\bar x\right\|_0}$ cells of dimension $s-\left\|\bar x\right\|_0$ attached, namely:
\[
 \bigcup_{\scriptsize \begin{array}{c}
        J \subset \left\{1, \ldots,n-\left\|\bar x\right\|_0\right\} \\
        1 \in J, |J| = s-\left\|\bar x\right\|_0 +1
   \end{array}} \mbox{conv} \left(e_j, j \in J\right).
\]
\end{theorem}

\proof
Theorem \ref{thm:def} allows deformations up to an arbitrarily small neighborhood of the M-stationary point $\bar x$. In such a neighborhood, we may assume without loss of generality that $\bar x=0$ and $f$ has the following form as from Theorem \ref{thm:morse}:
\begin{equation}
    \label{eq:n1}
        f(x)= f\left(\bar x\right) + \sum\limits_{i \in I_0\left(\bar x \right)}x_i + \sum\limits_{j \in I_1\left(\bar x \right)} x_j^2,
\end{equation}
where $x \in \R^{n-\left\|\bar x\right\|_0,s-\left\|\bar x\right\|_0}\times \R^{\left\|\bar x\right\|_0}$.

In terms of \cite{goresky:1988} the set 
$\R^{n-\left\|\bar x\right\|_0,s-\left\|\bar x\right\|_0}\times \R^{\left\|\bar x\right\|_0}$ can be interpreted as the product of the tangential part $\R^{\left\|\bar x\right\|_0}$ and the normal part $\R^{n-\left\|\bar x\right\|_0,s-\left\|\bar x\right\|_0}$. The cell-attachment along the tangential part is standard. Analogously to the unconstrained case, one cell of dimension zero has to be attached on $\R^{\left\|\bar x\right\|_0}$.  The cell-attachment along the normal part is more involved. Due to Lemma \ref{lem:cat}, we need to attach $\binom{n-\left\|\bar x\right\|_0-1}{s-\left\|\bar x\right\|_0}$ cells on $\R^{n-\left\|\bar x\right\|_0,s-\left\|\bar x\right\|_0}$, each of dimension $s-\left\|\bar x\right\|_0$. Finally, we apply Theorem 3.7 from Part I in \cite{goresky:1988}, which says that the local Morse data is the product of tangential and normal Morse data. Hence, the dimensions of the attached cells add together. This provides the assertion. \qed

Let us present a global interpretation of our results for SR. For that, we consider M-stationary points $\bar x$ with exactly $s-1$ non-zero entries, i.\,e.
\[
\left\|\bar x\right\|_0=s-1.
\]
We refer to them as saddle points.

\begin{theorem}[Morse relation for SR]
\label{thm:mrel}
Let $A$ be an $s$-regular matrix, and all M-stationary points of SR be nondegenerate with pairwise different functional values of the objective function. Then, it holds:
\begin{equation}
    \label{eq:mr-cs}
    (n-s)r_{1} \geq r-1,
\end{equation}
where $r$ is the number of local minimizers, and $r_1$ is the number of saddle points of SR.
\end{theorem}

\begin{figure}
    \centering
\begin{tikzpicture}
\draw[->] (0.5,0) -- (xyz cs:x=3);
\node at (2.8,-0.3){$x_3$};
\draw[->] (0,0.5) -- (xyz cs:y=3);
\node at (0.8,2.8){$x_{n-s+1}$};
\draw[->] (-0.35,-0.35) -- (xyz cs:z=5);
\node at (-2.0,-1.6){$x_1$};
\node at (0.8,0.8){$\ldots$};
\draw[] (-0.35,-0.35) -- (0.5,0);
\draw[] (-0.35,-0.35) -- (0,0.5);
\draw[densely dotted](0,0) -- (xyz cs:y=0.5);
\draw[densely dotted](0,0) -- (xyz cs:x=0.5);
\draw[densely dotted](0,0) -- (xyz cs:z=1);
\draw[densely dotted](0,0) -- (xyz cs:x=1, z=1);
\draw[->] (1,0,1) -- (4,0,4);
\node at (4,0.3,4){$x_2$};
\draw[->] (-0.35,-0.35) -- (1,0,1);
\fill (0.5,0) circle(2pt);
\fill (0,0.5) circle(2pt);
\fill (-0.35,-0.35) circle(2pt);
\fill (1,0,1) circle(2pt);
   \end{tikzpicture} 
   \caption{Cell-attachment in case 2)}
    \label{fig:2}
\end{figure}
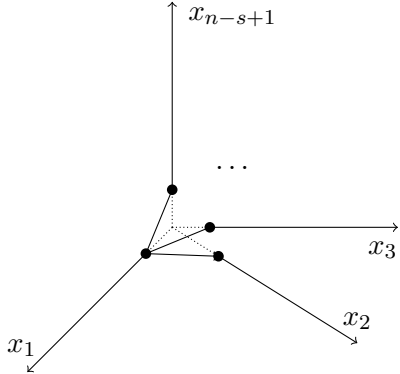

\proof 
Let $q_a$ denote the number of connected components of the lower level set $M^a$.
We focus on how $q_a$ changes as $a \in \R$ increases. Due to Theorem \ref{thm:def}, $q_a$ can change only if passing through a value corresponding to an M-stationary point $\bar x$, i.\,e. $a=f\left(\bar x\right)$. 
In fact, Theorem \ref{thm:def} allows homeomorphic deformations of lower level sets up to an arbitrarily small neighborhood of the M-stationary point $\bar x$.
Then, we have to estimate the difference between $q_a$ and $q_{a-\varepsilon}$, where $\varepsilon > 0$ is arbitrarily, but sufficiently small, and $a=f\left(\bar x\right)$. This is done by a {\it local argument}. 
We use Theorem \ref{thm:cell-a} which says that $M^{a}$ is homotopy-equivalent to $M^{a-\varepsilon}$ with a cell-attachment of
\begin{equation}
    \label{eq:hep4}
 \bigcup_{\scriptsize \begin{array}{c}
        J \subset \{1, \ldots,n-\left\|\bar x\right\|_0\} \\
        1 \in J, |J| = s-\left\|\bar x\right\|_0 +1
   \end{array}} \mbox{conv} \left(e_j, j \in J\right).
\end{equation}
Let us distinguish the following cases:
\begin{itemize}
    \item[1)] $\bar x$ is a local minimizer, i.\,e. $\left\|\bar x\right\|_0=s$.
    Then, by (\ref{eq:hep4}) we attach the cell $\mbox{conv}\left(e_1\right)$ of dimension zero to $M^{a-\varepsilon}$. Consequently, a new connected component is created, and it holds: 
\[
   q_a = q_{a-\varepsilon} + 1.
\]
    \item[2)] $\bar x$ is a saddle point, i.\,e. $\left\|\bar x\right\|_0=s-1$.
    Then, by (\ref{eq:hep4}) we attach $n-s+1$ cells of dimension one to $M^{a-\varepsilon}$, namely:
\[
   \bigcup_{\scriptsize \begin{array}{c}
        j=2, \ldots, n-s+1 
   \end{array}} \mbox{conv} \left(e_1,e_j\right).
\]
     Consequently, at most $n-s$ connected components disappear, and it holds:
\[
   q_{a-\varepsilon} - (n-s) \leq q_a \leq q_{a-\varepsilon}.
\]
For illustration we refer to Figure \ref{fig:2}.
    \item[3)] $\bar x$ is an M-stationary point with $\left\|\bar x\right\|_0<s-1$. The boundary of the cell-attachment in (\ref{eq:hep4}) is 
\[
  \bigcup_{\scriptsize \begin{array}{c}
        J \subset \{1, \ldots,n-\left\|\bar x\right\|_0\} \\
        1 \in J, |J| = s-\left\|\bar x\right\|_0+1 
   \end{array}} \partial \mbox{conv} \left(e_j, j \in J\right).
\]
The latter set is connected if $\left\|\bar x\right\|_0<s-1$. Consequently, the number of connected components of $M^{a}$ remains unchanged, and it holds:
\[
   q_a = q_{a-\varepsilon}.
\]
\end{itemize}

Now, we proceed with the {\it global argument}. Since 
the objective function is lower bounded by zero,
there exists $c<0$ such that $M^c$ is empty, thus, $q_c =0$. 
Due to Lemma \ref{lem:sr-fin}, the number of M-stationary points is finite. We conclude that there exists a lower level set $M^d$ which contains all, but finitely many, M-stationary points of SR. Due to Lemma \ref{lem:help}, it is possible to increase $d \in \R$ in order to additionally guarantee that $M^d$ is connected, i.\,e. $q_d=1$. Let us now vary the level $a$ from $c$ to $d$ and describe how the number $q_a$ of connected components of the lower level sets $M^a$ changes.
It follows from the local argument that $r$ new connected components are created, where $r$ is the number of local minimizers. Let $q$ denote the actual number of disappearing connected components if passing the levels corresponding to saddle points, and let $r_1$ denote the number of saddle points.
The local argument provides that at most $(n-s) r_{1}$ connected components might disappear while doing so, i.\,e.
\[
  q \leq (n-s) r_{1}.
\]
Altogether, we have:
\[
  r  - (n-s)r_{1} \leq r - q = q_d-q_c.
\]
By recalling that $q_d=1$ and $q_c=0$, we get Morse relation (\ref{eq:mr-cs}).
 \qed

We illustrate Theorem \ref{thm:mrel} by discussing the perturbed SR from Example \ref{ex:cw-inst}.

\begin{example}[Saddle point]
\label{ex:s-deg1}
Let the following sensing matrix and measurement vector be given:
\[
  A = \left( \begin{array}{cc}
    1   &  0\\
    0   &  1
   \end{array}\right), \quad b = \left( \begin{array}{c}
    1 \\ 1
   \end{array}\right).
\]
We consider the corresponding sparse  recovery problem with $s=1$:
\[
   \mbox{SR}: \quad \min_{x_1,x_2}\,\, \frac{1}{2}\left(x_1-1\right)^2 +\frac{1}{2}\left(x_2-1\right)^2 \quad \mbox{s.\,t.} \quad 
   \left\|\left(x_1, x_2\right)\right\|_0 \leq 1.
\]
As we have seen in Example \ref{ex:cw-inst}, both M-stationary points $(1,0)$ and $(0,1)$ are nondegenerate minimizers. Thus, we have $r=2$. Morse relation (\ref{eq:mr-cs}) from Theorem \ref{thm:mrel} provides:
\[
  r_1 \geq 1.
\]
Hence, there should exist a saddle point. In fact, $(0,0)$ is this nondegenerate saddle point, cf. Example \ref{ex:cw-inst}. Note that, due to $r_{1}=1$, Morse relation (\ref{eq:mr-cs}) holds with equality here.\qed
\end{example}

\section*{Acknowledgment}
The authors would like to thank Hubertus Th. Jongen for fruitful discussions.

\bibliographystyle{apalike}
\bibliography{lit.bib}

\end{document}